\title{A faithful tensor space representation \\ for the blob algebra} 
\author{P P Martin \\ \myaddress}\date{}
\documentclass[12pt]{article} 
\usepackage{epic} 
\usepackage{eepic} 
\usepackage{graphicx}
\usepackage{latexsym}
\usepackage{amssymb}
\providecommand{\noglossaryignore}[1]{}
\newcommand{\globalglossaryentry}[3]{\makebox[1.5in][l]{\tt $\backslash${#1}} 
\makebox[1.1in][l]{{$#2$}} \makebox[2.5in][l]{{#3}}\newline} 
\newcommand{\newcommandabbreviation}[3]{\newcommand{#1}{#2}%
\noglossaryignore{\globalglossaryentry{#3}{#2}{}}}
\newcommand{\renewcommandabbreviation}[3]{\renewcommand{#1}{#2}%
\noglossaryignore{\globalglossaryentry{#3}{#2}{}}}
\newcommand{\newcommandmacro}[4]{\newcommand{#1}{#2}%
\noglossaryignore{\globalglossaryentry{#3}{#2}{#4}}}
\newcommand{\gge}[3]{\noglossaryignore{\globalglossaryentry{#1}{#2}{#3}}}
\newcommand{\myaddress}%
{\parbox{3in}{\footnotesize \begin{center} 
Mathematics Department, City University, \\  
Northampton Square, London EC1V 0HB, UK.\end{center}}}
    
\newcounter{minidef}[section]

\newcounter{minicapt}

\newtheorem{de}{Definition}     \newtheorem{pr}{Proposition}

\noglossaryignore{GREEK ETC.\newline}
\newcommandabbreviation{\e}{\epsilon}{e}        
\newcommandabbreviation{\lam}{\lambda}{lam}  
\newcommandabbreviation{\la}{\langle}{la}        
\newcommandabbreviation{\ran}{\rangle}{ran}
\newcommandabbreviation{\ha}{\#}{ha}             
\newcommandabbreviation{\rmap}{\rightarrow}{rmap}
\newcommandabbreviation{\aaa}{\alpha}{aaa}        
\newcommandabbreviation{\ab}{\alpha,\beta}{ab}
\newcommandabbreviation{\aab}{a(\ab )}{aab}       
\noglossaryignore{\newline RINGS\newline}
\newcommandabbreviation{\HH}{H \!\!\! I}{HH}               
\newcommandabbreviation{\C}{\mathbb C}{C}
\newcommandabbreviation{\N}{\mathbb N}{N}   
\newcommandabbreviation{\Z}{\mathbb Z}{Z}      
\renewcommandabbreviation{\Re}{\mathbb R}{Re}
\newcommandabbreviation{\R}{{\mathbb R}}{R}
\newcommandabbreviation{\Q}{\mathbb Q }{Q}
\renewcommandabbreviation{\H}{\mathbb H }{H}
\noglossaryignore{\newline SYMMETRIC GROUP\newline}
\def\Sym(#1){\Sigma(#1)}                   
\gge{Sym(-)}{\Sym(-)}{symmetric group on - objects}
\def\Sy(#1){\Sigma_{#1}}                   
\gge{Sy(-)}{\Sy(-)}{symmetric group irreducible -}
\def\sym(#1){\mbox{\LARGE s}(#1)}        
\gge{sym(-)}{\sym(-)}{symmetric group on - objects (variant)}
\def\sy(#1){\mbox{\LARGE s}({#1})}        
\gge{sy(-)}{\sy(-)}{symmetric group irreducible - (variant)}
\newcommandmacro{\cs}{\C \, \sy(n)}{cs}{symmetric group algebra over $\C$}
\noglossaryignore{\newline PARTITIONS/SETS\newline}
\newcommand{\Nset}[1]{\underline{#1}}
\gge{Nset\{-\}}{\Nset{-}}{set of natural numbers to -} 
\def\nset(#1){ \{ #1 \}_{ \underline{n} }} 
\gge{nset(-)}{\nset(-)}{a set $-\times\Nset$}
\def\ul(#1){_{\underline{#1}}}             
\gge{ul(-)}{{}\ul(-)}{subscript underline -}
\def\Ee(#1){{\bf E}_{#1}}                  
\gge{Ee(-)}{\Ee(-)}{set of equivalence relations on set -}
\def\Eee(#1){{\bf E}_{\{ #1 \}_{\underline{n}}}}   
\gge{Eee(-)}{\Eee(-)}{ditto for nset}
\def\Een(#1,#2){{\bf E}_{\{ #1 \}_{\underline{#2}}}}   
\def\Ssn(#1,#2){{\bf S}_{\{ #1 \}_{\underline{#2}}}}   
\def\Ss(#1){{\bf S}_{#1}}                  
\def\Sss(#1){{\bf S}_{\{ #1 \}_{\underline{n}}}}   
\def\bbc(#1){((\beta_1)(\beta_2)...(\beta_{#1}))}      
\newcommandmacro{\Ln}{{\Gamma}^{n}}{Ln}{large index set}
\newcommandmacro{\LnQ}{{\Gamma}^{n}_Q}{LnQ}{index set}
\newcommandmacro{\Zz}{\zeta}{Zz}{`shape' function}
\noglossaryignore{\newline PARTITION ALGEBRA\newline}
\def\ka(#1){\kappa_{#1}}                   
\def\Sm(#1){\Sigma_{#1}}                   
\newcommandmacro{\com}{\bullet}{com}{bullet composition}
\newcommandmacro{\enm}{\; e^n(\! m\! ) \;}{enm}{product of idempotents}
\def\Ai(#1){ A^{ #1 \cdot } }              
\def\Aij(#1,#2){ A^{ #1  #2 } }            
\newcommandmacro{\One}{\mbox{\bf $1 \!\!\! 1$}}{One}{algebra unit 1}
\newcommandmacro{\Bp}{B_p}{Bp}{partition basis}
\def\Bb(#1){B_p[#1]}                       
\def\Pp(#1){P_n[#1]}                       
\def\Ps(#1){P_n[#1] \! /}                  
\newcommandmacro{\Ph}{\hat{P}}{Ph}{P hat  algebra}
\def\Is(#1){\sim^{#1}}                     
\noglossaryignore{\newline MODULES\newline}
\def\Wm(#1){{\cal S}_{#1}}                 
\gge{Wm(-)}{\Wm(-)}{Weyl module with index -}
\def\wm(#1,#2){{}_{#1}{\cal S}_{#2}}       
\gge{wm(-1,-)}{\wm(-1,-)}{Weyl module with index -}
\def\Ind(#1,#2,#3){\mbox{Ind}_{#1}^{#2}#3} 
\gge{Ind(-1,-2,-)}{\Ind(-1,-2,-)}{induction}
\def\Res(#1,#2,#3){\mbox{Res}_{#1}^{#2}#3} 
\gge{Res(-1,-2,-)}{\Res(-1,-2,-)}{restriction}
\newcommandabbreviation{\weyl}{standard}{weyl}
\newcommandabbreviation{\mod}{\mbox{mod}}{mod}
\newcommandabbreviation{\head}{\mbox{head }}{head}
\newcommandabbreviation{\Weyl}{Weyl}{Weyl}
\def\SS(#1){{\cal S}_{#1}}                 
\gge{SS(-)}{\SS(-)}{Specht/Weyl module index -}
\def\LL(#1){{\cal L}_{#1}}                 
\gge{LL(-)}{\LL(-)}{simple module index -}
\noglossaryignore{\newline FUNCTORS/MAPS\newline}
\newcommandmacro{\Gg}{{\cal G}}{Gg}{G Functor}
\newcommandmacro{\Fg}{{\cal F}}{Fg}{F Functor}
\newcommandmacro{\ra}{\rightarrow}{ra}{}
\def\ses(#1,#2,#3){0\ra #1 \ra #2 \ra #3 \ra 0}   
\gge{ses(1,2,-)}{\ses(1,2,-)}{\hspace{.5in} short exact sequence}
\def\starr(#1){ \stackrel{ #1 }{\longrightarrow} }
\gge{starr(-)}{\starr(-)}{}
\newcommandmacro{\doublerightarrow}{\; -\!\!\! -\!\!\!\!\!\! \gg \;}
{doublerightarrow}{}
\noglossaryignore{\newline PARTITION ALGEBRA MAPS\newline}
\newcommandmacro{\smap}{s}{smap}{`inclusion' map}
\newcommandmacro{\tmap}{t}{tmap}{$ P_n -> S_n$}
\newcommandmacro{\pmap}{\psi}{pmap}{$ S_n -> P_n $}
\noglossaryignore{\newline MISC.\newline}
\def\Amap(#1){{\cal A}_{#1}}               
\gge{Amap(-)}{\Amap(-)}{}
\def\Rr(#1){R_{#1}}                        
\gge{Rr(-)}{\Rr(-)}{restriction of E}
\def\Cr(#1){C_{#1}}                        
\gge{Cr(-)}{\Cr(-)}{restriction of E to N}
\newcommandmacro{\Tm}{{\cal T}}{Tm}{Transfer Matrix}
\def\On(#1){{\cal I}_{#1}}
\gge{On(-)}{\On(-)}{}
\newcommandmacro{\UU}{\underline{\sqcup}}{UU}{}  
\newcommandmacro{\UUU}{\sqcup}{UUU}{}  
\newcommandmacro{\Vq}{V_Q^{\otimes n}}{Vq}{Potts config. space}
\def\bs(#1,#2){\mbox{{\Large $\ast$}}^{#1}_{#2}}  
\gge{bs(-,-)}{\bs(-,-)}{general plumbing multiplier}
\newcommand{\ignore}[1]{}
\gge{ignore\{-\}}{\ignore{-}}{ignore argument!}
\def\choo(#1,#2){ \left( \begin{array}{c} #1 \\ #2 \end{array} \right) } 
\gge{choo(-1,-)}{\choo(-1,-)}{choose}
\newcommand{\Qed}{$\Box$}
\gge{Qed}{\mbox{\Qed}}{QED}
\def\staq(#1){\stackrel{#1}{=}}            
\gge{staq(-)}{\staq(-)}{}
\def\stam(#1){\stackrel{#1}{\rightarrow}}  
\gge{stam(-)}{\stam(-)}{}
\def\mat{ \left( \begin{array} }    
\def\tam{ \end{array}  \right) }
\gge{mat/tam}{...}{matrix delimiters}
\newcommand{\beq}{\begin{equation} }
\def\eql(#1){ \begin{equation} \label{#1} 
}
\newcommand{\eq}{\end{equation} }
\def\eqal(#1){\begin{eqnarray} \label{#1} }
\def\eqa{\end{eqnarray} }
\def\lab(#1){\label{#1}
}
\def\prl(#1){ \begin{pr} \label{#1} 
}
\def\del(#1){ \begin{de} \label{#1} 
}
\gge{smeq\{-\}}{...}{small equation}
\gge{fneq\{-\}}{...}{very small equation}
\noglossaryignore{\newline HECKE/BLOB\newline}
\newcommandmacro{\Hnq}{H_n(q)}{Hnq}{ * freestanding symbol}
\newcommandmacro{\Hn}{H_n}{Hn}{      *-mod etc.}
\newcommandmacro{\A}{{\cal A}}{A}{}
\newcommandmacro{\Cwts}{C}{Cwts}{}
\newcommandmacro{\CA}{{\cal A}}{CA}{}

\newcommandmacro{\calA}{{\cal A}}{calA}{}
\newcommandmacro{\modi}{\mbox{Mod} }{modi}{was mod not modi!}
\newcommandmacro{\Wgen}{{\Bbb S}}{Wgen}{}
\def\ol(#1){\overline{#1}}
\newcommandmacro{\st}{\mbox{St}}{st}{}
\def\CMult(#1,#2){(#1:#2)}
\def\CM(#1,#2){( #1 : #2 )}
\def\FMult#1,#2{(#1:#2)}
\def\CF#1,#2{(#1:#2)}

\newcommandmacro{\Top}{\mbox{Top}}{Top}{}
\newcommandmacro{\Soc}{\mbox{Soc}}{Soc}{}
\newcommandmacro{\Head}{\mbox{Head}}{Head}{}
\newcommandmacro{\Filt}{{\cal F}}{Filt}{}
\newcommandmacro{\Mod}{\mbox{mod}}{Mod}{}
\newcommandmacro{\Resi}{\mbox{Res }}{Resi}{was without i!}
\newcommandmacro{\Indi}{\mbox{Ind }}{Indi}{was without i!}
\def\RR(#1,#2){R^{#1}_{#2}}   
\def\TT(#1,#2){T^{#1}_{#2}}   

\newcommandmacro{\Ann}{\mbox{Ann}}{Ann}{}
\newcommandmacro{\Cen}{\mbox{Cen}}{Cen}{}
\newcommandmacro{\End}{\mbox{End}}{End}{}
\newcommandabbreviation{\semisimple}{semisimple}{semisimple}
\newcommandabbreviation{\Bratteli}{Bratteli}{Bratteli}
\newcommandabbreviation{\JBC}{Jones Basic Construction}{JBC}
\newcommandabbreviation{\pa}{partition algebra}{pa}
\newcommandabbreviation{\TM}{transfer matrix}{TM}
\newcommandabbreviation{\PM}{Potts model}{PM}
\newcommandabbreviation{\QSC}{quantum spin chain}{QSC}
\newcommandabbreviation{\Hamiltonian}{Hamiltonian}{Hamiltonian}
\newcommandabbreviation{\YS}{Young symmetrizer}{YS}

\newcommand{\beqa}{\begin{eqnarray}}%
\newcommand{\eeqa}{\end{eqnarray}}%

\newcommand{\U}{U}

\newcommand{\mU}{{\cal  U}}
\newcommand{\ym}{{m}}
\newcommand{\TL}{ Temperley--Lieb}
\newcommand{\chix}{\chi}
\newcommand{\egen}{e}

\begin{document} \maketitle 
 \newcommand{\ignoreifnotdraft}[1]{}
\ignoreifnotdraft{
\pagestyle{myheadings} \markboth{Draft}{\today}
}

\ignoreifnotdraft{ \newpage } 

\section{Introduction}
The blob algebra is a diagram algebra extending the Temperley--Lieb algebra in a
fairly natural way, which 
has a number of very nice 
properties (see \cite{CoxGrahamMartin03} for a review). 
Some time ago Martin and Woodcock \cite{MartinWoodcock01pre}
stumbled across a 
curious  
`tensor space' representation of the blob algebra, which turns out
\cite{MartinRyom02} to be a full tilting module \cite{Donkin93} 
in quasihereditary specialisations 
\cite{ClineParshallScott88,DlabRingel89b}. 
This raises the possibility of some intriguing new
developments in invariant theory (see \cite{MartinRyom02} for a
discussion). 
In the study of affine Hecke algebra representation theory 
it is also useful for technical reasons  
(see \cite{CoxGrahamMartin03,GrahamLehrer03}), 
to study the blob algebra, and the tensor space representation, in arbitrary
specialisations, including {\em non}-quasihereditary cases. 
In particular it is useful to know if the tensor space 
representation is {\em faithful} in arbitrary specialisations. 
In this paper we answer this question in the affirmative. 

We begin by assembling the machinery we will need in the more familiar
context of the Temperley--Lieb (TL) algebra. 
The ordinary tensor space representation here 
\cite{TemperleyLieb71,Baxter82,Jimbo85}
was shown to be faithful a long time ago
\cite{Martin92,DuParshallScott98}, 
and we use a similar method to \cite{Martin92} here. 
However we are able to implement it 
in such a way that it is applicable to representations subject only to
a relatively flexible set of conditions. 
Using this flexibility, 
we are eventually able to apply the
method to the blob algebra, thus obtaining a sufficient condition for blob
representations to be faithful. 

In the final section we recall the construction of the blob tensor
space representation, from which it is evident that it satisfies this
condition. 

\medskip


In this paper  $K$ is a ring, $x$ an invertible element in $K$, 
$q=x^2$,  and  $[n] = q^{n-1} + q^{n-3} + \ldots + q^{1-n}$. 
\newcommand{\Tp}{{T}}
\newcommand{\Td}{{\mathcal T}}
Define $\Tp_n^K$ to be the $K$--algebra with 
generators $\{ 1,\U_1,\ldots,\U_{n-1} \}$ and relations
\begin{eqnarray} 
 \U_i \U_i &=& (q+q^{-1}) \U_i     \label{TL001} \\
 \U_i \U_{i\pm 1} \U_i &=& \U_i    \label{TL002} \\
 \U_i \U_j &=& \U_j \U_i \hspace{1in} \mbox{($|i-j|\neq 1$)}  \label{TL003}
\end{eqnarray}


\section{Temperley--Lieb shenanigans}
\newcommand{\DD}{{\mathcal D}}%
\newcommand{\seq}{\mbox{seq}}%
\newcommand{\Dup}[1]{#1^{\cup}}%
\newcommand{\Ddown}[1]{#1_{\cap}}%
\newcommand{\TU}{{\mathcal U}}%
\newcommand{\braki}[1]{\langle #1 |}%
\newcommand{\keti}[1]{| #1 \rangle}%
\newcommand{\wa}{{\mathcal W}}%
For $n+m$ even, 
an $(n,m)$ TL diagram is a rectangular frame with $n$ nodes on the
northern and $m$ nodes on  the southern edge; 
the $n+m$ nodes are connected in pairs
by non--touching lines in the plane interior to the frame. 
Two such diagrams are identified if they partition the set of nodes
into pairs in the same way. 
The set of such diagram is denoted $\DD(n,m)$. 
Label the northern nodes $1,2,..,n$ and the southern nodes
$1',2',..,m'$. Say $(ij) \in D$ if nodes $i,j$ (primed, unprimed or
mixed) are connected in diagram $D$. 
Write $1 \in \DD(n,n)$ for the element such that $(ii') \in 1$ for all
$i$. 
Write $\TU_j \in \DD(n,n)$ for the element such that 
$(j \; j\! +\! 1), (j' \; (j\! +\! 1)') \in \TU_j$ and 
$(ii') \in \TU_j$ for all $i\neq j,j+1$. 
For example
\[
\TU_1  := \includegraphics{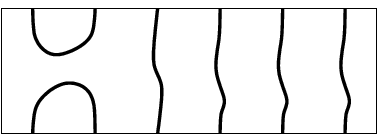}
\]
\newcommand{\dprod}[2]{#1  {\scriptsize \mbox{$\circ$}}  #2}
Define a product 
\begin{eqnarray*}
\DD(n,m) \times \DD(m,l) &\rightarrow& \DD(n,l)   \\
(D,D') &\mapsto& \dprod{D}{D'}
\end{eqnarray*}
by first concatenating the diagrams $D,D'$   
in such a way that the $i^{th}$ primed node
of $D$ meets the $i^{th}$ unprimed node of $D'$. 
(Call this object $D|D'$.) 
These nodes are then
discarded, leaving connections amongst the nodes of a resultant
diagram in $\DD(n,l)$. 
Note that $D|D'$ may have some
closed loops, which we ignore in $\dprod{D}{D'}$. However 
define a map 
\begin{eqnarray*}
\DD(n,m) \times \DD(m,l) &\rightarrow& \N   \\
(D,D') &\mapsto& \zeta((D,D'))
\end{eqnarray*}
where $\zeta((D,D'))$ is the number of closed loops discarded above. 
Thus for example $\dprod{\TU_1}{ \TU_1}=\TU_1$ and $\zeta((\TU_1,\TU_1))=1$.


The {\em propagating number} $\ha(D)$ of a diagram $D$ is the number
of lines of the form 
$(ij')$ in $D$. Note that it is possible to cut a diagram from the
western to the eastern edge in such a way that only these lines are
cut, and  they are cut once each. 
Let $\DD^l(n,m)$ denote the subset of $\DD(n,m)$ with propagating
number $l$. 
Note that cutting as above defines a unique map 
\begin{eqnarray*}
\DD^l(n,m) &\rightarrow& \DD(n,l)  \times \DD(l,m)  \\
D &\mapsto& (\Dup{D},\Ddown{D})
\end{eqnarray*}
such that $\dprod{\Dup{D}}{ \Ddown{D}} = D$. 


\newcommand{\deltap}{\delta'}
\newcommand{\qpo}[1]{\frac{\mbox{\tiny sign}(#1_i-#1_j)}{2}}
\newcommand{\cdepr}[4]{ \left( \!
    \prod_{ \scriptsize \begin{array}{c}i<j, \\ \! (#1)\in #4 \!\!\! \end{array}} 
    q^{\qpo{#2}}
    \delta'_{#2_i #3_j}
    \! \right)}%
\newcommand{\depr}[3]{ \left( 
    \prod_{                 
      \! (ij')\in #3 \!\!\! }
    \delta_{#1_i #2_j} \right)}%
\newcommand{\Rq}[3]{ \cdepr{ij}{#1}{#1}{#3} \!\!\!\! \depr{#1}{#2}{#3} \!\!\!\!
  \cdepr{i'j'}{#2}{#2}{#3} }%

Set $q=x^2$, $x \in K$, and define delta-function
\[
\delta_{ab} = \left\{ \begin{array}{ll} 1 & a=b \\ 
                                        0 & \mbox{otherwise}\end{array} \right.
\]
$\delta'_{ab}=1-\delta_{ab}$, and, for $a \in \Z \setminus \{0 \}$ 
\[
\mbox{sign}(a) = \left\{ \begin{array}{ll} +1 & a>0 \\ 
                                           -1 & a<0 \end{array}
                                       \right. .
\]
Associate to each $D\in \DD(n,m)$ a matrix $R_q(D)$ as follows. 
Rows are indexed by the set $\seq_n\{ 1,2 \}$ of words in $\{ 1,2 \}$ of length
$n$. Columns are indexed similarly by $\seq_m\{ 1,2 \}$. 
For $v \in \seq_n\{ 1,2 \}$ write $v_i$ for the $i^{th}$ term. 
Then 
\eql(R*)
R_q(D)_{vw} = \Rq{v}{w}{D}
\eq
For example, $R_q(1)$ is the unit matrix. 


\prl(agree)
Suppose there are no closed loops in $D|D'$.
Then for each pair $u,v$ there exists a unique $w$ giving rise 
to a non-vanishing summand in
$$(R_q(D)R_q(D'))_{uv} = \sum_w R_q(D)_{uw} R_q(D')_{wv} ;$$ and  
$$ R_q(D) R_q(D')= R_q( \dprod{D}{D'} )  . $$
More generally,
\eql(is a rep1)
 R_q(D) R_q(D') = [2]^{\zeta((D,D'))} R_q( \dprod{D}{D'} )  . 
\eq
\end{pr}
{\em Proof:} Fixing $u,v$ and
considering $\sum_w R_q(D)_{uw} R_q(D')_{wv}$ we have 
\[
\sum_w \Rq{u}{w}{D} 
\] 
\[ \hspace{2cm}
\Rq{w}{v}{D'}
\] 
Each delta-function factor corresponds to a  line in $D|D'$ (the
arguments correspond to the endpoints of the line). 
In particular each $w_i$ appears in two delta-functions. 
Hence each delta-function (or complementary delta-function) 
involving $w$ lies in a {\em chain} of one of a number of possible types. 
If there are no closed loops in $D|D'$ then 
those lines/deltas involving $w$ must lie in chains which begin either
in $u$ or in $v$. 
For example we might have $w_1,w_2$ appearing in the form 
$$\sum_{w_1, w_2} \delta_{u_1 w_1} q^{\qpo{w}} \delta'_{w_1 w_2} \delta_{w_2 u_2} 
=  q^{\qpo{u}} \delta'_{u_1 u_2}$$ 
where the right hand side shows the result of performing the relevant summations. 
Since every $w_i$ arises in this way, 
the complete sum may be replaced by precisely one term --- 
up to powers of $q$, a product of delta functions involving $u,v$. 
Considering an individual chain involving $w$,  
if it is ultimately propagating then an equal number of 
$(ij)$ lines from $D'$ and $(i'j')$ lines from $D$ 
are involved, contracting to a simple delta function. 
If it is ultimately within $u$ then there must be one more 
$(ij)$ line from $D'$ than $(i'j')$ lines from $D$, and so on.  

The general result follows
by a similar argument. 
\Qed


\begin{de} 
Two matrices $M,N$ are {\em mask equivalent} if $M_{ij}=0 \iff
N_{ij}=0$. Write $[M]$ for the equivalence class of $M$.  
\end{de}

Note, 
\eql(note*) [ R_q (D) ] = [R_{q'} (D)] . \eq
We will write $R(D)$ for $R_q(D)$ (other choices of parameter will be
written explicitly). 


\prl(claimX)
Provided there are no closed loops in $D|D'$, 
if $X \in [R(D)]$ and $Y \in [R(D')]$ then $XY \in [R( \dprod{D}{D'} )]$. 
\end{pr}
{\em Proof:} The delta function structure of $R(D)_{vw}$ has now been
overlain, in $X_{vw}$, with an arbitrary nonzero constant, 
$X_{vw} = k^{X}_{vw} R(D)_{vw}$, say. 
But since the delta function structure is the same, 
fixing $u,v$ we still have only one value of $w$ ($w^{{\tiny uv}}$ say) 
producing a nonvanishing term in
$\sum_w X_{uw} Y_{wv}$. 
Thus  
$\sum_w X_{uw} Y_{wv} 
= k^X_{u w^{{\tiny uv}}} k^Y_{w^{{\tiny uv}} v} \sum_w R(D)_{uw} R(D')_{wv}$. 
\Qed

\section{Temperley--Lieb algebra}

For $K$ a ring and $q$ a unit in $K$ let $\Td_n$ denote the \TL\ algebra,
a $K$--algebra 
with basis $\DD(n,n)$ and multiplication given by
\[
D.D' = [2]^{\zeta((D,D'))} \dprod{D}{D'} . 
\]
Thus from (\ref{is a rep1}), $R$ on $\DD(n,n)$ extends to a representation of
$\Td_n$ (in fact the usual action on tensor space 
\cite{TemperleyLieb71,Baxter82}). 
The following two results are standard
\cite{GoodmanDelaharpeJones89,Martin91}. 
\prl(genTL)
$\Td_n$ is generated by $\{ 1, \TU_1 , \TU_2 ,.., \TU_{n-1} \}$. 
\end{pr}
\prl(presenTL)
$\Td_n$ is isomorphic to the algebra with generators  
$\{ 1, U_1 , U_2 ,.., U_{n-1} \}$ and relations as in
equations~\ref{TL001} to~\ref{TL003}, 
with isomorphism given by $\TU_i \mapsto U_i$. 
\end{pr}


The Pascal triangle may be viewed as a graph embedded in the plane. 
It has vertices arranged in layers called levels. 
Levels are indexed $0,1,2,..$. 
Within level $i$ vertices are indexed by `column': $i,i-2,..,-i$. 
Thus a specific vertex may be labelled by (level,column)$=(i,i-2j)$. 
Edges are given by pairs of vertices: $((i,j),(i+1,j\pm 1))$. 
The 1--Pascal graph is the full subgraph on vertices with nonnegative
column index. 


Let $\wa_i(n)$ be the set of walks of length $n$ from (0,0) to
$(n,i)$ on the 1-Pascal graph. 
These walks may be represented in an obvious way by elements of
$\seq_n\{1,2\}$ (choose all such walks to start 1\ldots). 
Let $\wa^2_i(n) =  \wa_i(n) \times \wa_i(n)$ and 
$\wa^2(n) = \cup_i \wa_i(n) \times \wa_i(n)$. 
Draw an element
$(a,b)$ of $\wa^2(n)$ by drawing $a$ and the image of $b$ 
under reflection in the main
vertical of the Pascal triangle. 
The {\em envelope} of $(a,b) \in \wa^2_i(n)$ 
is the subset of the plane bounded by
this drawing and the piecewise straight 
line from vertex $(n,i)$ to $(n-i,0)$ to $(n,-i)$. 
For example, the envelope of $(121,112)$ is
\[
\includegraphics{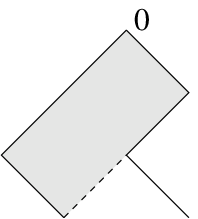}
\]

Partial order $\wa^2_i(n)$ by $(a,b) \leq (c,d)$ if 
the drawing of $(a,b)$ never leaves
the envelope formed by $(c,d)$. 
(We will also use the obvious underlying partial order on {\em single} walks.
This partial order is a lattice, with lowest walk 1212..1211..1, and
highest walk 11..122..2.) 
Extend to a 
partial order on $\wa^2(n)$ by $(a,b) \leq (c,d)$ if the endpoint 
of $a$ is $(n,i)$, that of $c$ is $(n,j)$, and $i<j$. 

The envelope of $(a,b)$ may be tiled in an obvious way with diamonds
(squares oriented at $45^o$) of side length 1. 
Form a map $$w: \wa^2(n) \rightarrow T_n$$ by scanning the tiling of
$(a,b)$ from top to bottom, left to right, and writing $U_i$ for each
tile with base at
row position $i$. 
Example:
\[
\includegraphics{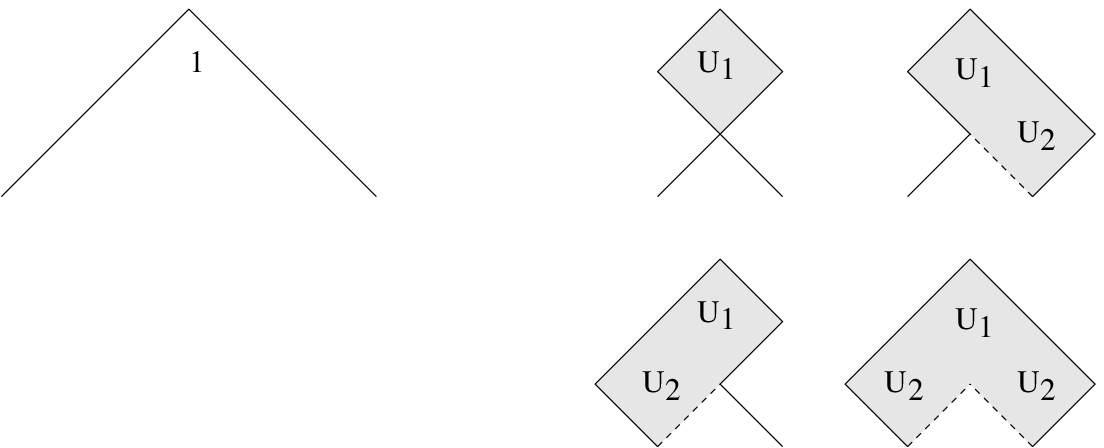}
\]

\prl(preTLmain1)
None of the elements  in $w( \wa^2(n))$ has a diagram representation
with a closed loop. As diagrams $w( \wa^2(n)) = \DD(n,n)$. 
\end{pr}
{\em Proof:} see for example \cite[\S 6.5]{Martin91}. 

\medskip


If $a$ is a walk or sequence with subsequence 21, 
with the 1 in the $i^{th}$ position, then write $a^i$ for the same
sequence except that the subsequence is replaces by 12. 
Note that $(a^i,b) > (a,b)$ for any $b$, 
and that $w((a^i,b))= U_i w((a,b))$. 


\prl(TLmain) (1) If $a=u$ and $b=v$ (confusing walks and sequences as
above) then
\[
R(w((a,b)))_{uv} \neq 0 . 
\]
(2) If 
\[
R(w((a,b)))_{uv} \neq 0 
\]
then $(a,b) \geq (u,v)$. 
\end{pr}
{\em Proof:} 
First note that 
(1) is true for the lowest walk pair in each lattice 
$\wa^2_i(n)$ by
an explicit calculation. For example, in bra--ket notation 
$$ \braki{1212..} U_1 U_3  \keti{1212..} = 
   \braki{1212..} U_3  \keti{q 1212.. + 2112..} =  $$ $$
   \braki{1212..}  \keti{q^2 1212.. + q 1221.. + q 2112.. + 2121..} =
   q^2 . $$
NB, (2) is the same as: if $(a,b)\not\geq(u,v)$ 
then $R(w((a,b)))_{uv} = 0$. 
Thus we may approach the whole proposition by working through various
   cases of $(a,b)$ and $(u,v)$. 
For our first case, 
suppose that $a$ ends at $(n,n-2i)$ and $u$ at $(n,n-2j)$ with $i>j$. 
In this case $(a,b)\not\geq(u,v)$ by virtue of their being 
in {\em different} lattices. 
Consider the lowest walk pair ($(a^o,a^o)$, say) in the lattice
   containing $(a,b)$. 
This has $w((a^o,a^o))=U_1 U_3..U_{2i-1}$. 
Given that $U_i \keti{..11..} = 0$, 
a simple sorting argument shows that there must be at least $i$ 2s in
the sequence $u$ for there to be a nonzero matrix element. 
In our case, however, there are precisely $j$ 2s. 
Thus $w((a^o,a^o))=U_1 U_3..U_{2i-1}$ vanishes on the whole
permutation block associated to $u$s of this type. 
But the rest of $w(\wa^2_{n-2i}(n))$ is
in the ideal generated by  $w((a^o,a^o))=U_1 U_3..U_{2i-1}$, so the image
of every pair in $\wa^2_{n-2i}(n)$ vanishes. 

It remains to deal with cases in which both $(a,b)$ and $(u,v)$ are
drawn from the same lattice $\wa^2_k(n)$ (some $k$).  
We work by induction on the lattice $w(\wa^2_k(n))$. 
That is, we suppose the proposition holds as regards 
all pairs below $(a,b)$, and all pairs $(u,v)$. 
Then in particular it holds 
for some pair $(a,c)$ such that $c^i=b$. We have 
\[
 \braki{a} w((a,b) ) \keti{b}
= \braki{a} w((a,c^i) ) \keti{c^i}
=\braki{a} w((a,c)) U_i \keti{c^i}
\] \[
=\braki{a} w((a,c))  \keti{q c^i +c}
=       q\braki{a} w((a,c))  \keti{ c^i }
        +\braki{a} w((a,c))  \keti{c}
\]
Since $(a,c) \not\geq (a,c^i)$ the first term vanishes 
by the inductive hypothesis; the second
does not, also by the inductive hypothesis. 
Thus  (1) holds provided the inductive step for (2) holds. 

As regards (2), first note that the base case is again straightforward: 
the lowest pair in
the lattice gives $U_1 U_3 ..$, which kills every sequence except 
the corresponding lowest one (the first step is always 1, $U_1$ kills
the sequence unless the second step is 2; the third step is now  forced to be 1,
and $U_3$ kills the sequence unless the fourth is 2; and so on). 
To prove the induction consider $\braki{u} w((a,b))$. 
By the left--right symmetry of our problem 
we are done if we can show this vanishes when 
$u \not\leq a$, so we restrict to such cases. 
We may assume WLOG that there is some $d$ and some $i$ such that
$a=d^i$, whereupon 
$$\braki{u} w((a,b)) = \braki{u} w((d^i,b)) = \braki{u} U_i w((d,b)) .
$$
(NB, $(d,b) < (a,b)$ so the inductive assumption holds for $(d,b)$ 
with all $(u,v)$.) 
If $u_i = u_{i+1}$ the last expression vanishes and we are done. 
Otherwise, we have 
$$\braki{u} U_i w((d,b))
= \left( q^{\pm 1} \braki{u} + \braki{u^{(i)}} \right) w((d,b))
$$
where  $u^{(i)}$ may be either higher or lower than $u$, depending on
whether $u=..12..$ or $..21..$ in the $i^{th}$ position. 
Since $u \not\leq a $ and $a> d$ we have $u \not\leq d $ and the first
term vanishes by the inductive assumption (note that the ket part is
not needed for this). 
If  $u^{(i)} > u$ then the second term vanishes similarly. 
If  $u^{(i)} < u$ then the $i^{th}$ and $i+1^{th}$ elements of both
$u^{(i)}$ and $d$ are 21. 
Thus $u \not\leq a $ implies $u^{(i)} \not\leq d $. 
\Qed

\prl(TLindep)
The matrices $R(w(\wa^2(n)))$ are a linearly independent set. 
The representation $R$ is faithful. 
\end{pr}
{\em Proof:} 
Pick a total order consistent with the partial order. 
As we run up through the order there is, for each
element, a matrix element which becomes nonzero first for that
element. 

It is a straightforward exercise to show that 
$|\wa^2(n) | = \mbox{Rank}(T_n)$. 
\Qed


Since the proof above uses only the occurences of nonzero matrix
elements we have, by the same token, a result on mask equivalent matrices:
\prl(maskver)
Any set $\{ X_D \in [R(D)] \; | \;  D \in w(\wa^2(n)) \}$ is linearly 
independent. 
\end{pr}
\section{The blob algebra}
\newcommand{\vargamma}{\gamma}%

In what follows it is convenient to 
shift the indices on the generators of $T_{2n}$ so that they run 
$U_{-n+1},U_{-n+2},..,U_{0},..,U_{n-1}$.
As before, the matrices $R(U_i)$ provide a representation of
$T_{2n}$. 


A line in a TL diagram is {\em exposed} if it may  be deformed to
touch the western edge of the frame. 
A {\em blob diagram} is like a TL diagram, except that any exposed line may
be decorated with a blob. For example
\[
e := \includegraphics{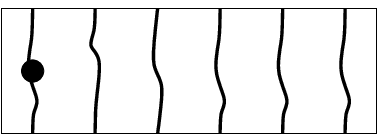}
\]
Write $\DD^b(n,m)$ for the extension of $ \DD(n,m)$ to include
decorated diagrams in this way. 
Blob diagram composition is  
like TL diagram composition, 
except that: 
\newline $\circ \;$ 
two blobs on the same line
may be replaced by one blob and a factor $\delta_e \in K$, and 
\newline $\circ \;$ 
a closed loop with a blob is replaced by a
factor $\vargamma \in K$ instead of $[2]$. 
\newline
Thus the blob algebra $b_n$ with basis $\DD^b(n,n)$ and 
this composition has three parameters, $q$, $\vargamma$ and $\delta_e$.
(Over a field, $\delta_e \neq 0$ may be rescaled to 1 without loss of
generality, but this point need not concern us here.)  

The proof of the following is straightforward. 
\prl(genblob)
$b_n$ is generated by $\{ 1, e, \TU_1 , \TU_2 , \ldots, \TU_{n-1} \}$. 
In particular, every element of $\DD^b(n,n)$ may be expressed as a
word in this set. 
\end{pr}


A word in the generators of $b_n$ (resp. $T_n$) is {\em loop free} if
its blob (resp. TL) diagram has no loops. 
This includes reduced words and words whose diagrams differ from those
of reduced words only by ambient isotopies. 
Let $B_n$ be a set of loop free words in the generators of $b_n$ which, 
as blob diagrams, form the diagram basis of $b_n$. 
Let $f$ be the map from words in the $b_n$ generators to words in the
$T_{2n}$ generators given 
piecewise 
by $f(e)=U_0$ and $f(\TU_i)=U_{-i}U_i$. 
Note that this takes loop free words to loop free words, 
but is not an algebra map. 
It induces an injective set map from the diagram basis of $b_n$ into
$\DD(2n,2n)$. 

\del(mirror)
A representation $\rho$ of $b_n$ 
with  $\rho(e) \in [R_r(U_0)]$, and $\rho(\TU_i)$ of the form $ X_i Y_i$, 
where $X_i \in [R_s(U_{-i})]$ and $Y_i \in [R_t(U_{i})]$ 
for some $r,s,t$, is called a {\em mirror} representation.
\end{de} 
(NB, by equation~\ref{note*} the choice of $r,s,t$ is actually
irrelevant to this statement; the point of introducing mask equivalence 
is that it does not differentiate $r,s,t$,
but preserves `enough' of the structure of diagram composition, as in
proposition~\ref{claimX}, to allow us to prove the following.) 
\prl(xxx)
Any mirror  representation of $b_n$ is faithful. 
\end{pr}
{\em Proof:}
Consider the set of matrices $\{ \rho(w) \; | \; w \in B_n \}$. 
It follows from proposition~\ref{claimX} that $\rho(w)$ 
is mask equivalent to $R(f(w))$.%
\footnote{The words $f(w)$ 
  correspond to TL diagrams which are left--right symmetric. 
  Indeed all symmetric diagrams in $\DD(2n,2n)$ may be obtained in this
  way (there are 
$ \frac{(2n)!}{n!n!}  
= |f(B_n)|$ of them).}
For example, $\rho(U_i) = X_i Y_i \in [ R(U_i U_{-i} )]$. 
At the level of diagrams, the map $f$ from $B_n$ 
to the set of TL diagrams is injective. 
That is, the set $f(B_n)$ of $T_{2n}$ words is, 
as a set of diagrams, a subset of $\DD(2n,2n)$. 
It therefore  follows from proposition~\ref{maskver} 
that the set of representation matrices for the diagram
basis of $b_n$ is linearly independent. 
\Qed
 



\section{Mirror representations}

We recall the representation $\rho_{0}$ of $b_n$ defined in
\cite[\S6.1]{MartinWoodcock01pre}: 
As explained in \cite{CoxGrahamMartin03}, the most interesting
unanswered questions about $b_n$ concern certain cases in which
$(\vargamma,\delta_e)$ 
can be written in the form $([m-1]\alpha,[m]\alpha)$ for some $m \in \Z$
and some scalar $\alpha$. 
(For example, $b_n$ is not quasihereditary in general over a field in
which $[2]=0$ and $\vargamma=0$ and $\delta_e =1$.)
Accordingly we recall the representation $\rho_{0}$ in an integral
form suitable for passing to such cases. 


Set
$$\mU^q(\chix)=\mat{cccc} 
0&0&0&0\\ 
0&q&1&0 \\
0&1&q^{-1}&0 \\
0&0&0&\chix \tam$$
and $\mU^q =\mU^q(0)$. 

Let $V_2 = K^2$. Fix $n$ and let $M_2^r(U_i) \in \End(V_2^{\otimes 2n})$ 
be a matrix acting trivially on every tensor
factor except the $i^{th}$ and $(i+1)^{th}$, where it acts as
$-\mU^r$. 
(Thus $M_2^r(U_i) = R_r(U_i)$ for $U_i \in \Tp_{2n}$ with $q=r$.)



Suppose there is an element $a \in K$ such that $a^4=-1$. 
Then $a^2 + a^{-2} = 0$. 
Fix $m$ such that $q^m \in K$ and 
set 
$$r= a^2 q^{\ym}$$ 
$$s  
 = a^5 x $$ 
$$t  
 = a^3 x $$ 

Let $b_n^{\Z[q,q^{-1}]}(q,\ym)$ be  
$\Z[q,q^{-1}] \DD^b(n,n) \subset b_n$
with $\vargamma = q^{m-1}-q^{-m+1}$ and $\delta_e = q^{m}-q^{-m}$, 
a $\Z[q,q^{-1}]$--algebra. 
Then there is an algebra homomorphism 
$$\rho_{0} : b_n^{\Z[q,q^{-1}]}(q,\ym) 
  \longrightarrow End_{\Z[a,x,x^{-1}]}(V_2^{\otimes 2n})$$ 
given by 
\beqa \label{mape}
\rho_{0}: \egen & \mapsto & 
                        a^{-2} M_2^r(U_n)
\\
\rho_{0}: U_i  & \mapsto &  M_2^s(U_{n-i}) M_2^t(U_{n+i}) . 
\eeqa
Comparing with definition~\ref{mirror} and (\ref{R*}) we see that $\rho_{0}$
is a mirror representation, and  hence faithful. 


\section{Discussion}
In \cite{Martin92} corresponding statements to proposition~\ref{TLindep} are
proved for each of the ordinary Hecke algebra quotients
$\End_{U_qsl_N}(V_N^{\otimes n})$, $V_N = K^N$ (explicitly for $K=\C$,
since this is a Physics reference, but the restriction is not
forced). It would be extremely desirable to generalise the blob
version in an analogous way, 
since the generalised blob algebras provide direct information about affine
Hecke representation theory \cite{MartinWoodcock01pre}. 
So far not even a candidate representation is known!


\bibliographystyle{amsplain}
\bibliography{new31,main,emma}

\end{document}